# Mean-field lattice trees

Christian Borgs and Jennifer Chayes
Microsoft Research
1 Microsoft Way
Redmond, WA 98052, USA
borgs@microsoft.com, jchayes@microsoft.com

Remco van der Hofstad[*] and Gordon Slade[†]
Department of Mathematics and Statistics
McMaster University
Hamilton, ON, Canada L8S 4K1
slade@math.mcmaster.ca

April 1, 1999

**Abstract**

We introduce a mean-field model of lattice trees based on embeddings into $\mathbb{Z}^d$ of abstract trees having a critical Poisson offspring distribution. This model provides a combinatorial interpretation for the self-consistent mean-field model introduced previously by Derbez and Slade, and provides an alternate approach to work of Aldous. The scaling limit of the mean-field model is integrated super-Brownian excursion (ISE), in all dimensions. We also introduce a model of weakly self-avoiding lattice trees, in which an embedded tree receives a penalty $e^{-\beta}$ for each self-intersection. The weakly self-avoiding lattice trees provide a natural interpolation between the mean-field model ($\beta = 0$), and the usual model of strictly self-avoiding lattice trees ($\beta = \infty$) which associates the uniform measure to the set of lattice trees of the same size.

## 1 Introduction

It is often the case that the understanding of a statistical mechanical model can be enhanced by the analysis of a corresponding *mean-field* model. The mean-field model typically involves an interaction that is simple enough to enable explicit computations of its scaling behaviour. Moreover, the scaling behaviour of the mean-field model is typically identical to that of the original model, above an upper critical dimension. The basic example is the Curie–Weiss model of ferromagnetism, which is the mean-field model corresponding to the Ising model [?]. For the self-avoiding walk, the mean-field model is simple random walk [?]. In both cases, the upper critical dimension is 4.

---
[*]Present address: TWI (SSOR), Delft University of Technology, Mekelweg 4, 2628 CD Delft, The Netherlands. Email: hofstad@twi.tudelft.nl

[†]Address after July 1, 1999: Department of Mathematics, University of British Columbia, Vancouver, BC, Canada V6T 1Z2. Email: slade@math.ubc.ca

In this paper, we introduce a mean-field model for lattice trees, based on lattice embeddings of abstract trees having a critical Poisson offspring distribution. The *scaling limit* of the mean-field model is the limit obtained by embedding increasingly large trees into an increasingly finer lattice. We will show that the scaling limit of the mean-field model is the random probability measure on $\mathbb{R}^d$ known as integrated super-Brownian excursion (ISE). This provides a simple construction of ISE. In this construction, we embed discrete abstract trees into $\mathbb{Z}^d$ and then take a continuum limit. Such an approach was outlined by Aldous in [?], although constructions more commonly first take a continuum limit of abstract trees and then embed these limiting objects into $\mathbb{R}^d$ [?, ?, ?]. Further discussion of ISE can be found in [?, ?, ?, ?, ?].

Our use of the Poisson distribution simplifies the analysis, but seems inessential. Although we do not prove convergence to ISE for other distributions, we will comment very briefly in Sections ?? and ?? on how our results might be generalised.

The occurrence of ISE as the scaling limit of lattice trees and of the incipient infinite percolation cluster above their respective upper critical dimensions 8 and 6 is discussed in [?, ?, ?, ?, ?, ?]. These interacting systems are more difficult to analyse than the mean-field model introduced here, and the methods described below serve as a basis for their analysis. Our primary goal in this paper is to isolate and describe these elementary methods, which underlie the work of [?, ?, ?, ?, ?, ?].

As we will indicate, our mean-field model provides a combinatorial interpretation for the self-consistent mean-field model introduced in [?]. Our mean-field model is also closely related to the mean-field model for lattice trees discussed in [?]. In [?, (5.15)–(5.17)], non-integer values are given for quantities that purportedly count embedded trees of various kinds. However, if we assume that the embeddings in [?] tacitly involve the Poisson weight factors we will introduce below, then we recover some of the results of [?].

It is often convenient to introduce a small parameter into a statistical mechanical model, as in the Domb–Joyce model of weakly self-avoiding walks [?]. We will define a model of weakly self-avoiding lattice trees, which associates a factor $e^{-\beta}$ to each self-intersection of an embedded tree. The weakly self-avoiding lattice trees interpolate between the mean-field model ($\beta = 0$) and the standard model of lattice trees in which all lattice trees of a given size and containing the origin are assigned equal probability ($\beta = \infty$).

We expect that the analysis of [?, ?, ?], which proves that the scaling limit is ISE for lattice trees in $\mathbb{Z}^d$ in dimensions $d \gg 8$ and for sufficiently "spread-out" lattice trees in dimensions $d > 8$, can be easily extended to obtain similar results for weakly self-avoiding lattice trees for $d > 8$ and $\beta \ll 1$. In fact, the weakly self-avoiding lattice trees should be easier to handle because the small parameter $\beta$ is more explicit than the small parameters in [?, ?, ?]. However, we have not carried out the exercise of checking that the various calculations involved all go through in this setting.

## 2 Results

A *bond* in $\mathbb{Z}^d$ is a pair $\{x, y\}$ of sites $x, y \in \mathbb{Z}^d$ with $\|x - y\|_1 = 1$. A *lattice tree* is a finite connected set of bonds that contains no cycles. We will say that $x$ is in a lattice tree $L$ if there is a bond in $L$ that contains $x$. We put the uniform measure on the set of all $n$-site lattice trees containing the origin. This model is difficult to analyse because of the self-avoidance constraint inherent in the prohibition on cycles. The mean-field model will relax this restriction completely.

The mean-field model is defined in terms of embeddings of abstract trees into $\mathbb{Z}^d$. The abstract trees are the family trees of the critical birth process with Poisson offspring distribution. In more detail, we begin with a single individual having $\xi$ offspring, where $\xi$ is a Poisson random variable of mean 1, i.e., $\mathbb{P}(\xi = m) = (em!)^{-1}$. Each of the offspring then independently has offspring of its own, with the same critical Poisson distribution. For a tree $T$ consisting of exactly $n$ individuals, with the $i^{\text{th}}$ individual having $\xi_i$ offspring, this associates to $T$ the weight

$$\mathbb{P}(T) = \prod_{i \in T} \frac{e^{-1}}{\xi_i!} = e^{-n} \prod_{i \in T} \frac{1}{\xi_i!}. \tag{2.1}$$

The product is over the vertices of $T$.

It is important to be clear about when two trees $T$ are the same and when they are not. For this, we introduce a description of $T$ in terms of *words*. These words arise inductively as follows. The root is the word $0$. The children of the root are the words $01, 02, \ldots 0\xi_0$. The children of $01$ are the words $011, \ldots, 01\xi_{01}$, and so on. The family tree is then uniquely represented by a set of words. Two trees are the same if and only if they are represented by the same set of words. In the terminology of [?], we are dealing with *plane trees*.

We define an embedding $\varphi$ of $T$ into $\mathbb{Z}^d$ to be a mapping from the vertices of $T$ into $\mathbb{Z}^d$, such that the root is mapped to the origin and adjacent vertices in the tree are mapped to nearest neighbours in $\mathbb{Z}^d$. Given a tree $T$ having $|T|$ vertices, there are $(2d)^{|T|-1}$ possible embeddings $\varphi$ of $T$. The mean-field model is then defined to be the set of configurations $(T, \varphi)$, with probabilities

$$\mathbb{P}(T, \varphi) = \frac{1}{(2d)^{|T|-1}} \mathbb{P}(T). \tag{2.2}$$

Equivalently, we may regard the mean-field model as corresponding to branching random walk, with mean-1 Poisson branching distribution.

The one-point and two-point functions of the mean-field model are defined for complex $z$ and $\zeta$, with $|z|, |\zeta| \leq 1$, by

$$t_z^{(1)} = \sum_{(T,\varphi)} \mathbb{P}(T,\varphi) z^{|T|} = \sum_T \mathbb{P}(T) z^{|T|}, \tag{2.3}$$

$$t_{z,\zeta}^{(2)}(x) = \sum_{(T,\varphi)} \mathbb{P}(T,\varphi) z^{|T|} \sum_{i \in T} I[\varphi(i) = x] \zeta^{|i|} \quad (x \in \mathbb{Z}^d), \tag{2.4}$$

where $|i|$ denotes the graph distance from $i$ to the root of $T$. The series (??) converges for $|z| \leq 1$, with $t_1^{(1)} = 1$. The series (??) clearly converges for $|z| < 1$, $|\zeta| \leq 1$. The one-point function is a generating function for embedded trees rooted at the origin, while the two-point function is a generating function for embedded trees rooted at the origin and containing the site $x$.

The one-point function is given by the following theorem.

**Theorem 2.1** *For $d \geq 1$, the one-point function is given by $t_z^{(1)} = \sum_{n=1}^{\infty} \frac{n^{n-1}}{n!} e^{-n} z^n$, which is one solution of the implicit equation*

$$t_z^{(1)} e^{-t_z^{(1)}} = z e^{-1}. \tag{2.5}$$

Equation (??) actually defines a function analytic in $\mathbb{C} \setminus [1, \infty)$, and we are taking the principal branch. Equation (??) can be written in terms of the Lambert $W$ function, defined by $We^W = z$, as $t_z^{(1)} = -W(-ze^{-1})$. The branches of $W$ are described in [?]. Theorem ?? rederives the well-known result that for critical Poisson branching processes,

$$\mathbb{P}(|T| = n) = \frac{n^{n-1}}{n!} e^{-n}. \tag{2.6}$$

Given a summable function $f : \mathbb{Z}^d \to \mathbb{C}$, its Fourier transform is given by $\hat{f}(k) = \sum_{x \in \mathbb{Z}^d} f(x) e^{ik \cdot x}$. We write

$$\hat{D}(k) = \frac{1}{d} \sum_{j=1}^{d} \cos k_j, \quad k = (k_1, \ldots, k_d) \in [-\pi, \pi]^d \tag{2.7}$$

for the Fourier transform of the step distribution for the simple random walk on $\mathbb{Z}^d$ taking nearest-neighbour steps with equal probabilities. The two-point function of the mean-field model is then given by the following theorem.

**Theorem 2.2** *For $d \geq 1$, $k \in [-\pi, \pi]^d$, $|z| < 1$, $|\zeta| \leq 1$,*

$$\hat{t}_{z,\zeta}^{(2)}(k) = \frac{t_z^{(1)}}{1 - t_z^{(1)} \zeta \hat{D}(k)}. \tag{2.8}$$

*The denominator of the right side vanishes for $z = \zeta = 1$, $k = 0$, and in that case $\hat{t}_{1,1}^{(2)}(0) = \infty$.*

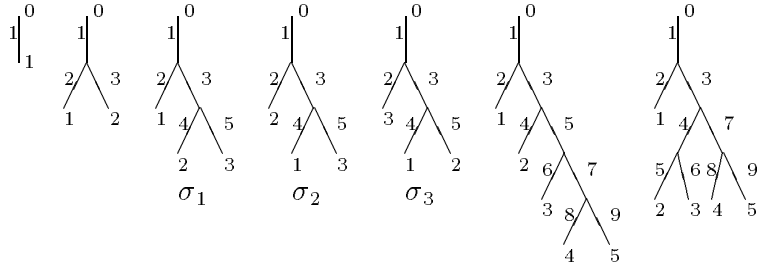

Figure 1: The shapes for $m = 2, 3, 4$, and examples of the $7!! = 7 \cdot 5 \cdot 3 = 105$ shapes for $m = 6$. The shapes' edge labellings are arbitrary but fixed.

Apart from unimportant factors, the one-point and two-point functions given above are the same as those of the mean-field models of [?] and [?]. The two-point function given in Theorem ?? can be interpreted as the two-point function of simple random walk with an activity $\zeta$ associated to each step of the walk and an activity $t_z^{(1)}$ associated to each site. We may therefore regard an embedded tree containing 0 and $x$ as corresponding to a simple random walk path from 0 to $x$ with a one-point function attached at each site along the way.

In order to define $m$-point functions, for $m \geq 3$, we first introduce the notion of *shape*. We start with an abstract $m$-skeleton, which is a tree having $m$ unlabelled external vertices of degree 1 and $m - 2$ unlabelled internal vertices of degree 3, and no other vertices. An $m$-*shape* is a tree having $m$ labelled external vertices of degree 1 and $m - 2$ unlabelled internal vertices of degree 3, and no other vertices, *i.e.*, an $m$-shape is a labelling of an $m$-skeleton's external vertices by the labels $0, 1, \ldots, m - 1$. When $m$ is clear from the context, we will refer to an $m$-shape simply as a shape. For notational convenience, we associate to each $m$-shape an arbitrary labelling of its $2m - 3$ edges, with labels $1, \ldots, 2m - 3$. This arbitrary choice of edge labelling is fixed once and for all. Thus an $m$-shape $\sigma$ is a labelling of an $m$-skeleton's external vertices together with a corresponding specification of edge labels. Let $\Sigma_m$ denote the set of $m$-shapes. There is a unique shape for $m = 2$ and $m = 3$, and $(2m - 5)!!$ distinct shapes for $m \geq 4$ (see [?, (5.96)] for a proof). In this notation, $(-1)!! = 1$ and $(2j + 1)!! = (2j + 1)(2j - 1)!!$ for $j \geq 0$.

Next, we need the notion of backbone. We write $\bar{\imath} = (i_1, \ldots, i_{m-1})$ for a sequence of $m - 1$ vertices $i_j$ in a tree $T$ (possibly with repetition), and define the *backbone* $B$ of $(T, \bar{\imath})$ to be the subtree of $T$ spanning $0, i_1, \ldots, i_{m-1}$. There is an induced labelling of the external vertices of the backbone, in which vertex $i_l$ is labelled $l$. Ignoring vertices of degree 2 in $B$, this backbone is equivalent to a shape $\sigma_B$ or to its modification by contraction of one or more edges to a point. (In the latter case, as we will discuss further in Section ??, the choice of $\sigma_B$ may not be unique.) The edge labels of $\sigma_B$ induce labels on the paths in $T$ comprising the backbone $B$.

Finally, we need a notion of compatibility. Let $\vec{s} = (s_1, \ldots, s_{2m-3})$ for nonnegative integers $s_j$, and $\vec{y} = (y_1, \ldots, y_{2m-3})$ for $y_j \in \mathbb{Z}^d$. Note that we distinguish $m - 1$ and $2m - 3$ component vectors by using $\bar{\phantom{x}}$ and $\vec{\phantom{x}}$ respectively. Restoring vertices of degree 2 in $B$, let $b_j$ denote the length of the backbone path corresponding to edge $j$ of $\sigma_B$, with $b_j = 0$ for any contracted edge. We say that $(T, \varphi, \bar{\imath})$ is *compatible* with $(\sigma; \vec{y}, \vec{s})$ if $\sigma_B$ can be chosen such that $\sigma_B = \sigma$, if $b_j = s_j$ for all edges $j$ of $\sigma$, and if the image under $\varphi$ of the backbone path (oriented away from the root) corresponding to $j$ undergoes the displacement $y_j$ for all edges $j$ of $\sigma$. Note that $\vec{y}$ determines the image under $\varphi$ of all $2m - 2$ vertices of $\sigma_B$.

Now we define the $m$-point function (with $m - 2$ additional internal vertices) by

$$t_n^{(m)}(\sigma; \vec{y}, \vec{s}) = \sum_{(T,\varphi):|T|=n} \mathbb{P}(T, \varphi) \sum_{i_1,\ldots,i_{m-1} \in T} I[(T, \varphi, \bar{\imath}) \text{ is compatible with } (\sigma; \vec{y}, \vec{s})]. \quad (2.9)$$

We also define

$$t_{z,\vec{\zeta}}^{(m)}(\sigma; \vec{y}) = \sum_{n=0}^{\infty} \sum_{s_1,\ldots,s_{2m-3}=0}^{\infty} t_n^{(m)}(\sigma; \vec{y}, \vec{s}) z^n \prod_{j=1}^{2m-3} \zeta_j^{s_j}, \quad (2.10)$$

and the Fourier transform

$$\hat{f}(\vec{k}) = \sum_{y_1,\ldots,y_{2m-3} \in \mathbb{Z}^d} f(\vec{y}) e^{i\vec{k} \cdot \vec{y}}, \quad (2.11)$$

where $\vec{k} \cdot \vec{y} = \sum_{j=1}^{2m-3} k_j \cdot y_j$.

**Theorem 2.3** *For $d \geq 1$, $m \geq 2$, $k_j \in [-\pi, \pi]^d$, $|z| < 1$, $|\zeta_j| \leq 1$,*

$$\hat{t}_{z,\vec{\zeta}}^{(m)}(\sigma; \vec{k}) = \left(t_z^{(1)}\right)^{-2(m-2)} \prod_{j=1}^{2m-3} \hat{t}_{z,\zeta_j}^{(2)}(k_j) = t_z^{(1)} \prod_{j=1}^{2m-3} \frac{1}{1 - t_z^{(1)} \zeta_j \hat{D}(k_j)}. \quad (2.12)$$

Theorem ?? gives the same $m$-point function as the self-consistent mean-field model of [?], apart from the factor $\left(t_z^{(1)}\right)^{-2(m-2)}$ that was absent in [?]. This factor has a natural combinatorial interpretation. Namely, it "corrects" for an overcounting of the branch at each of the $m - 2$ internal shape vertices, as this branch is counted in $\prod_{j=1}^{2m-3} \hat{t}_{z,\zeta_j}^{(2)}(k_j)$ once by each of the three two-point functions incident at that vertex. This factor is equal to 1 at the critical point $z = 1$, and does not play a role in the leading scaling behaviour.

We now turn to the scaling limit of the mean-field model. A discussion of the scaling limit requires a digression concerning super-processes, a topic that has received considerable attention in probability theory [?, ?]. The most basic example of a super-process is super-Brownian motion, which arises as the scaling limit of branching random walk. Super-Brownian motion is a continuous-time stochastic process in which the state space is the set of finite measures on $\mathbb{R}^d$. It dies out in finite time, and its total mass integrated over its entire history is a random variable. Integrated super-Brownian excursion (ISE) is the random probability measure on $\mathbb{R}^d$ obtained by conditioning this total mass, integrated over time, to be 1. We will consider the scaling limit of mean-field lattice trees conditioned to have mass $n$. This is the same as conditioning on the mass of branching random walk, and it is to be expected that the scaling limit will therefore be ISE. We will prove that this is indeed the case. In fact, as we will discuss in Section ??, our proof provides an elementary construction of ISE.

To state our result, we write $M_1(\dot{\mathbb{R}}^d)$ for the compact space of probability measures on the one-point compactification $\dot{\mathbb{R}}^d$ of $\mathbb{R}^d$, equipped with the topology of weak convergence (see [?]). We regard the set $M_1(\mathbb{R}^d)$ of probability measures on $\mathbb{R}^d$ as embedded in $M_1(\dot{\mathbb{R}}^d)$. ISE is a certain probability measure $\mu_{\text{ISE}}$ on $M_1(\mathbb{R}^d)$, i.e., it is the law of a random probability measure on $\mathbb{R}^d$. Mean-field lattice trees arising from embeddings of trees $T$ with $|T| = n$ induce a measure $\mu_n$ on $M_1(\mathbb{R}^d)$, as follows. Given a tree $T$ with $|T| = n$, and an embedding $\varphi$ of $T$, let $\mu(T, \varphi)$ be the probability measure on $\mathbb{R}^d$ which assigns mass $n^{-1} \sum_{i \in T} I[\varphi(i) = x]$ ($x \in \mathbb{Z}^d$) to points $xd^{1/2}n^{-1/4}$ in $\mathbb{R}^d$. The measure $\mu_n$ on $M_1(\mathbb{R}^d)$ is then the measure that assigns mass $\mathbb{P}((T, \varphi)| |T| = n)$ to each $\mu(T, \varphi)$ with $|T| = n$. In other words, we obtain a random probability measure on $\mathbb{R}^d$ by assigning equal mass to each of the $n$ embedded vertices of a rescaled version of an embedded tree.

**Theorem 2.4** *For $d \geq 1$, as measures on $M_1(\dot{\mathbb{R}}^d)$, $\mu_n$ converges weakly to $\mu_{\text{ISE}}$.*

The weak convergence in Theorem ?? is the assertion that for any continuous function $F$ on $M_1(\dot{\mathbb{R}}^d)$, $\lim_{n \to \infty} \int_{M_1(\dot{\mathbb{R}}^d)} F(\nu) d\mu_n(\nu) = \int_{M_1(\dot{\mathbb{R}}^d)} F(\nu) d\mu_{\text{ISE}}(\nu)$. A result along the lines of Theorem ?? was already sketched in [?].

Next, we introduce a model of weakly self-avoiding lattice trees. For $\beta \geq 0$, let

$$Z_n^\beta = \sum_{(T,\varphi): |T|=n} \mathbb{P}(T, \varphi) \exp\left[-\tfrac{1}{2}\beta \sum_{i,j \in T: i \neq j} I[\varphi(i) = \varphi(j)]\right], \quad (2.13)$$

and, for $|T| = n$, define

$$\mathbb{Q}_n^\beta(T, \varphi) = \frac{1}{Z_n^\beta} \mathbb{P}(T, \varphi) \exp\left[-\tfrac{1}{2}\beta \sum_{i,j \in T: i \neq j} I[\varphi(i) = \varphi(j)]\right]. \quad (2.14)$$

The measure $\mathbb{Q}_n^\beta$ on the set of embedded $n$-site trees rewards self-avoidance by giving a penalty $e^{-\beta}$ to each self-intersection of an embedded tree. For $\beta = 0$, $\mathbb{Q}_n^0$ is just our mean-field model conditional on $|T| = n$. The next theorem shows that the weakly self-avoiding lattice trees interpolate between the mean-field model and lattice trees, in the sense that $\mathbb{Q}_n^\infty$ corresponds to the uniform measure on the set of $n$-site lattice trees containing the origin. In the statement of the theorem, $\ell_n$ denotes the number of $n$-site lattice trees containing the origin. Given an injective $\varphi$ and a lattice tree $L$, we abuse notation by writing $\varphi(T) = L$ if $\varphi(T)$ consists of the sites in $L$ and the edges in $T$ are mapped to the bonds in $L$.

**Theorem 2.5** *For $d \geq 1$ and $n \geq 0$, $\lim_{\beta \to \infty} \mathbb{Q}_n^\beta(T, \varphi) = 0$ if $\varphi$ is not injective. Given an $n$-site lattice tree $L$, $\lim_{\beta \to \infty} \sum_{(T,\varphi): \varphi(T)=L} \mathbb{Q}_n^\beta(T, \varphi) = \ell_n^{-1}$.*

The remainder of this paper is organised as follows. In Section ??, we prove Theorems ??, ?? and ??. In Section ??, we prove Theorem ?? and provide an additional statement concerning convergence of backbones. Finally, in Section ??, we prove Theorem ??.

# 3 The $m$-point functions

In this section, we prove Theorems ??, ?? and ??. To indicate the role of the Poisson distribution, we begin the proofs of Theorems ?? and ?? with a general critical offspring distribution $p_m = \mathbb{P}(\xi = m)$ and then specialize to the critical Poisson distribution $p_m = 1/(em!)$. For the general offspring distribution, (??) becomes

$$\mathbb{P}(T,\varphi) = \frac{1}{(2d)^{|T|-1}}\mathbb{P}(T) = \frac{1}{(2d)^{|T|-1}} \prod_{i \in T} p_{\xi_i}. \tag{3.1}$$

We introduce the generating function $P(w) = \sum_{m=0}^{\infty} p_m w^m$. Note that $P(1) = 1$, and since the offspring distribution is assumed to be critical, $P'(1) = 1$. For the critical Poisson distribution, $P(w) = e^{w-1}$. We will make use of the fact that the Poisson distribution has moments of all orders, but we expect that a more careful analysis can be used to prove that mean-field lattice trees with a critical offspring distribution having finite variance will converge to ISE in the scaling limit. A general result of this form is stated in [?]. It would be of interest to prove this using our methods, but we will prove Theorem ?? in Section ?? only for the Poisson distribution.

For the proofs, it will be helpful to introduce generating functions for planted plane trees (in the terminology of [?]). A planted plane tree is a plane tree for which the root has exactly one offspring. Generating functions for planted plane trees will be denoted by $r$ rather than $t$, and, by convention, will not include the factor $zp_1$ associated with the root having exactly one offspring. Thus the one-point function for planted plane trees is given simply by

$$r_z^{(1)} = t_z^{(1)}. \tag{3.2}$$

The one-point function on the right side of (??) arises as the generating function of the root's child and its progeny.

**Proof of Theorem ??.** Conditioning on the number of offspring of the root gives

$$t_z^{(1)} = \sum_{m=0}^{\infty} p_m z \left(r_z^{(1)}\right)^m = zP(t_z^{(1)}). \tag{3.3}$$

For the Poisson distribution, this implies $t_z^{(1)} = ze^{t_z^{(1)}-1}$, which gives (??). The Taylor expansion then follows from Lagrange's inversion formula (see, $e.g.$, [?, p.23]). □

We define the two-point function $r_{z,\zeta}^{(2)}(x)$ for planted plane trees to be the restriction of the summation in (??) to trees $T$ for which the root has a single offspring, with the factor $ze^{-1}$ associated with the root omitted. Then

$$\hat{r}_{z,\zeta}^{(2)}(k) = \zeta \hat{D}(k) \hat{t}_{z,\zeta}^{(2)}(k). \tag{3.4}$$

**Proof of Theorem ??.** The Fourier transform of the two-point function is given by

$$\hat{t}_{z,\zeta}^{(2)}(k) = \sum_{(T,\varphi)} \sum_{j \in T} \mathbb{P}(T,\varphi) z^{|T|} e^{ik \cdot \varphi(j)} \zeta^{|j|}. \tag{3.5}$$

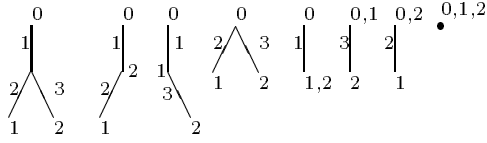

Figure 2: The $2^3 = 8$ subshapes for $m = 3$.

The contribution to the right side arising when $j$ is the root is simply $t_z^{(1)}$. When $j$ is not the root, we condition on the number of offspring of the root. With (??), this gives

$$\hat{t}_{z,\zeta}^{(2)}(k) = t_z^{(1)} + \sum_{m=1}^{\infty} p_m z m \left(r_z^{(1)}\right)^{m-1} \hat{r}_{z,\zeta}^{(2)}(k) = t_z^{(1)} + zP'(t_z^{(1)})\hat{r}_{z,\zeta}^{(2)}(k). \qquad (3.6)$$

In the middle of (??), the factor $z$ takes care of the factor associated with the root on the left side, and the factor $m$ corresponds to choosing which of the root's offspring is an ancestor of the vertex $j$. By (??), this gives

$$\hat{t}_{z,\zeta}^{(2)}(k) = \frac{t_z^{(1)}}{1 - zP'(t_z^{(1)})\zeta\hat{D}(k)}. \qquad (3.7)$$

Specializing now to the Poisson distribution, we have $zP'(t_z^{(1)}) = zP(t_z^{(1)}) = t_z^{(1)}$, which gives the desired result. Note that $\hat{t}_{1,1}^{(2)}(0) = \sum_{(T,\varphi)} \sum_{i \in T} \mathbb{P}(T, \varphi) = \sum_n n\mathbb{P}(|T| = n) = \infty$, by (??) and Stirling's formula. □

For $m \geq 3$, we introduce the notion of a *subshape* of a shape $\sigma \in \Sigma_m$. A subshape of $\sigma \in \Sigma_m$ is an abstract tree obtained by contracting a subset of the edges of $\sigma$ to a point. This can lead to multiply-labelled vertices, and contracted edges lose their labelling. The subshapes for $m = 3$ are shown in Figure ??. In general, there are $2^{2m-3}$ subshapes of a shape $\sigma \in \Sigma_m$. We denote subshapes by $\lambda$ and write $\lambda \leq \sigma$ when $\lambda$ is a subshape of $\sigma$. We denote the set of edge labels of $\lambda$ by $e(\lambda)$.

**Proof of Theorem ??.** For the Poisson offspring distribution, (??) gives $\hat{t}_{z,\zeta}^{(2)}(k) = t_z^{(1)}[1 + \hat{r}_{z,\zeta}^{(2)}(k)]$. Therefore, it suffices to show that

$$\hat{t}_{z,\vec{\zeta}}^{(m)}(\sigma; \vec{k}) = t_z^{(1)} \prod_{j=1}^{2m-3} \left(1 + \hat{r}_{z,\zeta_j}^{(2)}(k_j)\right). \qquad (3.8)$$

Expanding the product, the desired identity (??) is equivalent to

$$\hat{t}_{z,\vec{\zeta}}^{(m)}(\sigma; \vec{k}) = t_z^{(1)} \sum_{\lambda \leq \sigma} \prod_{j \in e(\lambda)} \hat{r}_{z,\zeta_j}^{(2)}(k_j). \qquad (3.9)$$

Given a subshape $\lambda$, we let $t(\lambda)$ denote the result of restricting the summation in (??) to $s_j = 0$ for $j \notin e(\lambda)$ and $s_j > 0$ for $j \in e(\lambda)$. Its Fourier transform will be denoted $\hat{t}(\lambda)$.

We leave implicit the dependence on the variables of $\vec{k}$ and $\vec{\zeta}$, as these are determined by the edge labels of $\lambda$. Then

$$\hat{t}_{z,\vec{\zeta}}^{(m)}(\sigma; \vec{k}) = \sum_{\lambda \leq \sigma} \hat{t}(\lambda). \tag{3.10}$$

Thus it suffices to show that

$$\hat{t}(\lambda) = t_z^{(1)} \prod_{j \in e(\lambda)} \hat{r}_{z,\zeta_j}^{(2)}(k_j). \tag{3.11}$$

We will use $\pi$ to denote a subshape for which the root has exactly one offspring, and write $\hat{r}(\pi) = (zp_1)^{-1}\hat{t}(\pi)$. As before, the factor $(zp_1)^{-1}$ serves to cancel the factor $zp_1$ associated to the root in $\hat{t}(\pi)$. We denote by $\bar{\pi}$ the subshape obtained from $\pi$ by contracting the edge incident on the root. We claim that

$$\hat{r}(\pi) = \hat{r}_{z,\zeta}^{(2)}(k) \frac{1}{t_z^{(1)}} \hat{t}(\bar{\pi}), \tag{3.12}$$

where $\zeta$ and $k$ bear the subscript of the label of the edge incident on the root of $\pi$. The identity (**??**) will be proved below. Given a subshape $\lambda$ having at least one edge, let $\pi_1, \ldots, \pi_b$ be the branches emerging from its root. For a general offspring distribution, as in (**??**),

$$\hat{t}(\lambda) = \sum_{j=b}^{\infty} zp_m j(j-1)\cdots(j-b+1)\left(r_z^{(1)}\right)^{j-b} \prod_{a=1}^{b} \hat{r}(\pi_a) = zP^{(b)}(t_z^{(1)}) \prod_{a=1}^{b} \hat{r}(\pi_a). \tag{3.13}$$

For the Poisson distribution, $zP^{(b)}(t_z^{(1)}) = t_z^{(1)}$, independently of $b$. This leads to the simple result $\hat{t}(\lambda) = t_z^{(1)} \prod_{a=1}^{b} \hat{r}(\pi_a)$, and the desired result (**??**) then follows by substituting (**??**) into this identity recursively. For general offspring distributions, degree-dependent vertex factors will arise in (**??**).

It remains to prove (**??**). We do this by conditioning on whether the length of the tree's backbone path, corresponding to the edge of $\pi$ incident on the root, is equal to or greater than 1. This leads, by conditioning as in (**??**), to

$$\hat{r}(\pi) = \zeta \hat{D}(k)\hat{t}(\bar{\pi}) + \zeta \hat{D}(k) t_z^{(1)} \hat{r}(\pi). \tag{3.14}$$

Solving and using (**??**) and (**??**), we obtain

$$\hat{r}(\pi) = \frac{\zeta \hat{D}(k)}{1 - \zeta \hat{D}(k) t_z^{(1)}} \hat{t}(\bar{\pi}) = \hat{r}_{z,\zeta}^{(2)}(k) \frac{1}{t_z^{(1)}} \hat{t}(\bar{\pi}). \tag{3.15}$$

$\square$

# 4 Convergence to ISE

In Sections ?? and ??, we prove Theorem ??. Section ?? contains an additional statement concerning convergence of backbones. As we will explain at the end of Section ??, our proof of Theorem ?? provides a construction of the measure $\mu_{\text{ISE}}$.

## 4.1 Moment measures

To prove weak convergence of probability measures on $M_1(\dot{\mathbb{R}}^d)$, it is sufficient to prove weak convergence of the moment measures [?, Lemma 2.4.1(b)]. The $l^{\text{th}}$ moment measure $M^{(l)}$ of a probability measure $\mu$ on $M_1(\dot{\mathbb{R}}^d)$ is the probability measure on $(\dot{\mathbb{R}}^d)^l$ defined, for $l \geq 1$, by $dM^{(l)}(x_1, \ldots, x_l) = \int_{M_1(\dot{\mathbb{R}}^d)} d\mu(\nu) d\nu(x_1) \cdots d\nu(x_l)$. To prove weak convergence of the moment measures, it is sufficient to prove pointwise convergence of their characteristic functions. We begin by introducing the ISE moment measures.

Let $m \geq 2$. Given a shape $\sigma \in \Sigma_m$, we associate to edge $j$ (oriented away from vertex 0) a nonnegative real number $t_j$ and a vector $y_j$ in $\mathbb{R}^d$. Writing $\vec{y} = (y_1, \ldots, y_{2m-3})$ and $\vec{t} = (t_1, \ldots, t_{2m-3})$, we define

$$a^{(m)}(\sigma; \vec{y}, \vec{t}) = \left(\sum_{j=1}^{2m-3} t_j\right) e^{-(\sum_{j=1}^{2m-3} t_j)^2/2} \prod_{j=1}^{2m-3} \frac{1}{(2\pi t_j)^{d/2}} e^{-y_j^2/2t_j} \tag{4.1}$$

and

$$A^{(m)}(\sigma; \vec{y}) = \int_0^\infty dt_1 \cdots \int_0^\infty dt_{2m-3} \, a^{(m)}(\sigma; \vec{y}, \vec{t}). \tag{4.2}$$

Then $\int_{\mathbb{R}^{d(2m-3)}} A^{(m)}(\sigma; \vec{y}) d\vec{y} = 1/(2m-5)!!$, so the sum of this integral over shapes $\sigma \in \Sigma_m$ is equal to 1. Let $\vec{k} \cdot \vec{y} = \sum_{j=1}^{2m-3} k_j \cdot y_j$, with each $k_j \in \mathbb{R}^d$. The Fourier integral transform $\hat{A}^{(m)}(\sigma; \vec{k}) = \int_{\mathbb{R}^{d(2m-3)}} A^{(m)}(\sigma; \vec{y}) e^{i\vec{k} \cdot \vec{y}} d\vec{y}$ is given by

$$\hat{A}^{(m)}(\sigma; \vec{k}) = \int_0^\infty dt_1 \cdots \int_0^\infty dt_{2m-3} \, \hat{a}^{(m)}(\sigma; \vec{k}, \vec{t}), \tag{4.3}$$

with

$$\hat{a}^{(m)}(\sigma; \vec{k}, \vec{t}) = \left(\sum_{j=1}^{2m-3} t_j\right) e^{-(\sum_{j=1}^{2m-3} t_j)^2/2} \prod_{j=1}^{2m-3} e^{-k_j^2 t_j/2}. \tag{4.4}$$

The functions (??) and (??) are further discussed in [?] (see also [?, ?, ?, ?, ?]).

The $l^{\text{th}}$ moment measure $M^{(l)}$ for ISE can be written in terms of $A^{(l+1)}$, for $l \geq 1$. This is a deterministic measure which is absolutely continuous with respect to Lebesgue measure on $\mathbb{R}^{dl}$. The first moment measure $M^{(1)}$ has density $A^{(2)}(x)$. The second moment measure $M^{(2)}$ has density $\int A^{(3)}(y, x_1 - y, x_2 - y) d^d y$. In general, the density of $M^{(l)}$ at $x_1, \ldots, x_l$, for $l \geq 3$, is given by integrating $A^{(l+1)}(\sigma, \vec{y})$ over $\mathbb{R}^{d(l-1)}$ and then summing over the $(2l-3)!!$ shapes $\sigma$. Here $\vec{y}$ consists of integration variables $y_j$ corresponding to the edges $j$ on paths from vertex 0 to vertices of degree 3 in $\sigma$, and the other $y_a$ are fixed by the requirement that each external vertex $x_i$ is given by the sum of the $y_e$ over the edges $e$ connecting vertices 0 and $i$ in $\sigma$. Thus, the integration corresponds to integrating over the $l-1$ internal vertices, with the $l+1$ external vertices fixed at $0, x_1, \ldots, x_l$. For example, the contribution to the density of $M^{(3)}$ due to $\sigma_1$ of Figure **??** is $\int A^{(4)}(\sigma_1; y_1, x_1 - y_1, y_3, x_2 - y_1 - y_3, x_3 - y_1 - y_3) d^d y_1 d^d y_3$.

The characteristic function of $M^{(l)}$ can be written in terms of the functions (**??**). For $l = 1$, $\hat{M}^{(1)}(k) = \hat{A}^{(2)}(k)$. Similarly, for $l = 2$, there is a single shape and $\hat{M}^{(2)}(k_1, k_2) = \int A^{(3)}(y, x_1 - y, x_2 - y) e^{ik_1 \cdot x_1} e^{ik_2 \cdot x_2} d^d y d^d x_1 d^d x_2 = \hat{A}^{(3)}(k_1 + k_2, k_1, k_2)$. For $l \geq 3$, there is more than one shape, and

$$\hat{M}^{(l)}(\bar{k}) = \sum_{\sigma \in \Sigma_{l+1}} \hat{A}^{(l+1)}(\sigma; \vec{k}) \tag{4.5}$$

with each of the $2m - 3$ components of $\vec{k}$ given by a specific linear combination (depending on $\sigma$) of the $l$ components of $\bar{k} = (k_1, \ldots, k_l)$. For example, for $l = 3$ and the shape $\sigma_1$ of Figure **??**, $(\sigma_1; \vec{k}) = (\sigma_1; k_1 + k_2 + k_3, k_1, k_2 + k_3, k_2, k_3)$.

Next, we introduce the moment measures $M_n^{(l)}$ of the mean-field lattice trees. For $l \geq 1$, let

$$s_n^{(l+1)}(x_1, \ldots, x_l) = \sum_{(T, \varphi): |T| = n} \mathbb{P}(T, \varphi) \sum_{i_1, \ldots, i_l \in T} \prod_{j=1}^{l} I[\varphi(i_j) = x_j]. \tag{4.6}$$

Recall the definition of $\mu_n$ above Theorem **??**. Writing $\bar{x} = (x_1, \ldots, x_l)$, the $l^{\text{th}}$ moment measure $M_n^{(l)}$ of $\mu_n$ is the deterministic probability measure on $\mathbb{R}^{dl}$ which places mass $[n^l \mathbb{P}(|T| = n)]^{-1} s_n^{(l+1)}(\bar{x})$ at $\bar{x} d^{1/2} n^{-1/4}$, for $\bar{x} \in \mathbb{Z}^{dl}$. The characteristic function $\hat{M}_n^{(l)}(k)$ of $M_n^{(l)}$ is given by

$$\hat{M}_n^{(l)}(\bar{k}) = \frac{1}{n^l \mathbb{P}(|T| = n)} \hat{s}_n^{(l+1)}(\bar{k} d^{1/2} n^{-1/4}), \tag{4.7}$$

where $\bar{k} \cdot \bar{x} = k_1 \cdot x_1 + \cdots + k_l \cdot x_l$ and

$$\hat{f}(\bar{k}) = \sum_{\bar{x}} f(\bar{x}) e^{i\bar{k} \cdot \bar{x}}. \tag{4.8}$$

To prove Theorem **??**, it suffices to prove that

$$\lim_{n \to \infty} \hat{M}_n^{(l)}(\bar{k}) = \hat{M}^{(l)}(\bar{k}), \quad (l \geq 1). \tag{4.9}$$

In fact, proving (**??**) provides a construction of ISE, as we now explain. Because $M_1(\dot{\mathbb{R}}^d)$ is compact, some subsequence of the sequence $\mu_n$ converges to a limit $\mu$. Given (**??**), it follows that $\mu$ has moments $M^{(l)}$, which uniquely characterises the measure $\mu$. But it then follows, again from (**??**), that $\mu_n$ must converge to $\mu$. Thus it follows from (**??**) that a limiting measure $\mu$ exists, and this provides a construction of $\mu_{\text{ISE}} = \mu$. Initially, this constructs $\mu_{\text{ISE}}$ as a measure on $M_1(\dot{\mathbb{R}}^d)$. However, since the moment measures $M^{(l)}$ have no mass at the point at infinity, $\mu_{\text{ISE}}$ is in fact a measure on $M_1(\mathbb{R}^d)$.

## 4.2 Convergence of moment measures

In this section, we complete the proof of Theorem **??** by proving (**??**). Define $t_n^{(m)}(\sigma;\vec{y})$, for $m \geq 2$, by

$$t_{z,\vec{1}}^{(m)}(\sigma;\vec{y}) = \sum_{n=1}^{\infty} t_n^{(m)}(\sigma;\vec{y})z^n. \tag{4.10}$$

For $m = 2, 3$, there is only one shape and, suppressing $\sigma$ in the notation, we have

$$s_n^{(2)}(x) = t_n^{(2)}(x), \quad s_n^{(3)}(x_1, x_2) = \sum_y t_n^{(3)}(y, x_1 - y, x_2 - y). \tag{4.11}$$

The relation between $s_n^{(m)}$ and $t_n^{(m)}$, for $m \geq 4$, will be discussed below. The basic ingredient of the proof of (**??**) is the following lemma. In its statement, the notation $f(n) \sim g(n)$ means that $\lim_{n \to \infty} f(n)/g(n) = 1$.

**Lemma 4.1** *For $m \geq 2$, $k_j \in \mathbb{R}^d$, and $n \to \infty$,*

$$\hat{t}_n^{(m)}(\sigma; \vec{k}d^{1/2}n^{-1/4}) \sim \frac{1}{\sqrt{2\pi}} n^{m-5/2} \hat{A}^{(m)}(\sigma; \vec{k}). \tag{4.12}$$

**Proof.** By Cauchy's theorem,

$$\hat{t}_n^{(m)}(\sigma; \vec{k}) = \frac{1}{2\pi i} \oint \hat{t}_{z,\vec{1}}^{(m)}(\sigma; \vec{k}) \frac{dz}{z^{n+1}}, \tag{4.13}$$

where the integral is around a circle of radius less than 1, centred at the origin. By Theorem **??**,

$$\hat{t}_{z,\vec{1}}^{(m)}(\sigma; \vec{k}) = t_z^{(1)} \prod_{j=1}^{2m-3} \frac{1}{1 - t_z^{(1)} \hat{D}(k_j)}. \tag{4.14}$$

By Taylor's theorem,

$$\hat{D}(k) = 1 - \frac{k^2}{2d} + O(k^4), \tag{4.15}$$

as $k \to 0$. Also, using (**??**) it can be shown that

$$t_z^{(1)} = 1 + O(|1-z|^{1/2}) \quad \text{and} \quad t_z^{(1)} = 1 - \sqrt{2(1-z)} + O(|1-z|), \tag{4.16}$$

with the error terms uniform in $\mathbb{C}\setminus[1, \infty)$. Substituting (**??**)–(**??**) into (**??**) gives

$$\hat{t}_{z,\vec{1}}^{(m)}(\sigma; \vec{k}d^{1/2}n^{-1/4}) = [1 + O(|1-z|^{1/2})] \prod_{j=1}^{2m-3} \frac{2}{n^{-1/2}k_j^2 + 2^{3/2}\sqrt{1-z} + O(n^{-1}k_j^4 + |1-z|)}. \tag{4.17}$$

It is then an exercise in contour integration, as in [**?**, Section 4.2], to deform the contour in (**??**) to the branch cut $[1, \infty)$ of the square root to conclude (**??**). □

By Stirling's formula and (??),

$$\mathbb{P}(|T|=n) = (2\pi)^{-1/2} n^{-3/2}[1+O(n^{-1})]. \tag{4.18}$$

Convergence of the first and second moments, i.e., (??) for $l=1,2$, then follows immediately from (??) and Lemma ??. The higher moments require further discussion.

Before considering the third and higher moments, we note that the above $n^{-3/2}$ behaviour is associated more generally with the size distribution of Galton–Watson trees whose offspring distribution has finite variance (see [?, Proposition 24]). This behaviour is associated with the $\sqrt{1-z}$ appearing in (??) and is thus closely connected with the occurrence of ISE as the scaling limit. This is consistent with the statement in [?] that the scaling limit is ISE for more general offspring distributions having finite variance. However, we are treating only the Poisson case here.

We will now consider the higher moments. For $l \geq 3$, Lemma ?? implies that

$$\lim_{n\to\infty} \frac{\sum_{\sigma \in \Sigma_{l+1}} \hat{t}_n^{(l+1)}(\sigma; \vec{k}d^{1/2}n^{-1/4})}{n^l \mathbb{P}(|T|=n)} = \sum_{\sigma \in \Sigma_{l+1}} \hat{A}^{(l+1)}(\sigma; \vec{k}). \tag{4.19}$$

If it were the case that $\hat{s}_n^{(l+1)}(\bar{k})$ were equal to $\sum_{\sigma \in \Sigma_{l+1}} \hat{t}_n^{(l+1)}(\sigma; \vec{k})$, convergence of all moments would be immediate. But $\hat{s}_n^{(l+1)}(\bar{k})$ is not equal to $\sum_{\sigma \in \Sigma_{l+1}} \hat{t}_n^{(l+1)}(\sigma; \vec{k})$, because it is not the case that $s_n^{(m)}(\bar{x})$ is equal to the sum of $t_n^{(m)}(\sigma; \vec{y})$ over all $(\sigma; \vec{y})$ that are compatible with $\bar{x}$ in the sense that the $x_i$ are given by the sum of the $y_j$ as prescribed by the shape $\sigma$. The discrepancy arises from degenerate trees whose backbone corresponds to a strict subshape with at least one shape edge contracted.

For example, there is a unique tree $T$ having just two vertices, i.e., the tree in which the root has one child and there are no further descendants. Thus $s_2^{(4)}(0,0,e_1) = (2d)^{-1}e^{-2}$, where $e_1 = (1,0,\ldots,0)$. However, this mean-field lattice tree containing the sites $x_1 = x_2 = 0$, $x_3 = e_1$, contributes $(2d)^{-1}e^{-2}$ to each of three choices of $t_2^{(4)}(\sigma; \vec{y})$, namely to $t_2^{(4)}(\sigma_1; 0,0,0,0,e_1)$, $t_2^{(4)}(\sigma_2; 0,0,0,0,e_1)$, and $t_2^{(4)}(\sigma_3; 0,e_1,0,0,0)$; see Figure ??. Thus it is not the case, in general, that $s_n^{(l+1)}(\bar{x})$ is given by the sum of $t_n^{(l+1)}(\sigma; \vec{y})$ over all $(\sigma; \vec{y})$ compatible with $\bar{x}$.

In view of (??), (??), (??) and (??), to prove convergence of the $l^{\text{th}}$ moments, for $l \geq 3$, it suffices to show that

$$\left| \hat{s}_n^{(l+1)}(\bar{k}) - \sum_{\sigma \in \Sigma_{l+1}} \hat{t}_n^{(l+1)}(\sigma; \vec{k}) \right| \leq O(n^{l-2}), \tag{4.20}$$

where $\vec{k}$ is determined by $\bar{k}$ and $\sigma$ as in (??). This difference then constitutes an error term, down by $n^{-1/2}$ compared to $\hat{s}_n^{(l+1)}(\bar{k})$, by Lemma ??. The remainder of the proof is devoted to obtaining (??).

Let $l \geq 3$, and suppose $T \ni i_1,\ldots,i_l$. If the backbone of $(T,\bar{\imath})$ corresponds to a full $(l+1)$-skeleton with no contracted edges, then $\bar{\imath}$ determines a labelling of an $(l+1)$-skeleton and therefore uniquely determines an $(l+1)$-shape. Whether or not the backbone corresponds to a full $(l+1)$-skeleton, given a shape consistent with the backbone, the $2l-1$ displacements $\vec{y}$ (possibly zero) consistent with that shape are uniquely determined by an embedding of $T$. Degeneracy of $(\sigma;\vec{y})$ thus requires the backbone to correspond to a strict subshape, and in that case, the maximum possible number of compatible choices for $(\sigma;\vec{y})$ is the number of shapes, which is $(2l-3)!!$.

We now write $s_n^{(l+1)}(\bar{x}) = u_n^{(l+1)}(\bar{x}) + e_n^{(l+1)}(\bar{x})$, where $u_n$ comprises the contributions to $s_n$ from full $(l+1)$-skeletons and $e_n$ comprises the contributions to $s_n$ from degenerate skeletons. It follows from the above discussion that, for $l \geq 3$,

$$\left| \hat{s}_n^{(l+1)}(\vec{k}) - \sum_{\sigma \in \Sigma_{l+1}} \hat{t}_n^{(l+1)}(\sigma;\vec{k}) \right| \leq [(2l-3)!! - 1]\hat{e}_n^{(l+1)}(\bar{0}). \tag{4.21}$$

It suffices to argue that the right side of (??) is at most $O(n^{l-2})$. For this, we introduce the generating function $E^{(l+1)}(z) = \sum_n \hat{e}_n^{(l+1)}(\bar{0})z^n$. This is a sum of terms of the form $\hat{t}(\lambda)$, where $\lambda$ is a strict subshape and all $k_j = 0$, $\zeta_j = 1$. By (??), (??) and (??), it follows that $|E^{(l+1)}(z)| \leq O(|1-z|^{-(l-1)})$ uniformly in $|z| < 1$, where the power $l-1 = \frac{1}{2}(2l-2)$ arises because at least one of the $2l-1$ backbone paths is trivial. Then [?, Lemma 3.2(i)] or [?, Theorem 4] implies the desired bound $\hat{e}_n^{(l+1)}(\bar{0}) \leq O(n^{l-2})$.

This completes the proof of Theorem ??.

## 4.3 Backbone convergence

The following lemma can be proved exactly as in [?, Section 4.2], so we omit the proof here.

**Lemma 4.2** *For $m \geq 2$, $\vec{u} = (u_1, \ldots, u_{2m-3})$ with $u_j \in [0, \infty)$, $\vec{k} = (k_1, \ldots, k_{2m-3})$ with $k_j \in \mathbb{R}^d$, and for $n \to \infty$,*

$$\hat{t}_n^{(m)}(\sigma; \vec{k}d^{1/2}n^{-1/4}, \lfloor \vec{u}n^{1/2} \rfloor) \sim \frac{1}{\sqrt{2\pi}} \frac{1}{n} \hat{a}^{(m)}(\sigma; \vec{k}, \vec{u}). \tag{4.22}$$

In (??), the backbone scaling is $n^{1/2}$. Since space is being scaled as $n^{1/4}$, this is Brownian scaling. The lemma provides an interpretation of the variables $u_j$ in $\hat{a}^{(m)}(\sigma; \vec{k}, \vec{u})$ as rescaled Brownian time variables along backbone paths.

# 5 Weakly self-avoiding lattice trees

**Proof of Theorem ??.** The first statement of the theorem, for non-injective $\varphi$, follows immediately from the definition of $\mathbb{Q}_n^\beta$.

For the second statement of the theorem, let $\mathcal{L}_n$ denote the set of $n$-site lattice trees containing the origin. This has cardinality $\ell_n$. We will prove that

$$\sum_{(T,\varphi): \varphi(T)=L} \mathbb{P}(T,\varphi) = (2d)^{-(n-1)} e^{-n} \tag{5.1}$$

for every $L \in \mathcal{L}_n$. The important point for our proof is that the right side is the same for all $L \in \mathcal{L}_n$, and its precise value plays no role. In fact, given (**??**), we then have

$$Z_n^\infty = \sum_{L \in \mathcal{L}_n} \sum_{(T,\varphi):\varphi(T)=L} \mathbb{P}(T,\varphi) = \ell_n (2d)^{-(n-1)} e^{-n}, \tag{5.2}$$

which gives the desired result that

$$\sum_{(T,\varphi):\varphi(T)=L} \mathbb{Q}_n^\infty(T,\varphi) = \frac{1}{Z_n^\infty} \sum_{(T,\varphi):\varphi(T)=L} \mathbb{P}(T,\varphi) = \frac{1}{\ell_n}. \tag{5.3}$$

To prove (**??**), we first note that by (**??**) and (**??**),

$$\sum_{(T,\varphi):\varphi(T)=L} \mathbb{P}(T,\varphi) = (2d)^{-(n-1)} e^{-n} \sum_{(T,\varphi):\varphi(T)=L} \prod_{i \in T} \frac{1}{\xi_i!}, \tag{5.4}$$

where $\xi_i$ is the number of offspring of vertex $i$. It suffices to show that

$$\sum_{(T,\varphi):\varphi(T)=L} \prod_{i \in T} \frac{1}{\xi_i!} = 1. \tag{5.5}$$

Let $b_0$ be the degree of $0$ in $L$, and given nonzero $x \in L$, let $b_x$ be the degree of $x$ in $L$ minus 1. Then the set $\{b_x : x \in L\}$ must be equal to the set of $\xi_i$ for any $T$ that can be mapped to $L$. Defining $\nu(L) = \#\{(T,\varphi) : \varphi(T) = L\}$, (**??**) is therefore equivalent to

$$\nu(L) = \prod_{x \in L} b_x!. \tag{5.6}$$

We prove (**??**) by induction on the number $N$ of generations of $L$. By this, we mean the length of the longest self-avoiding path in $L$, starting from the origin. The identity (**??**) clearly holds if $N = 0$. Our induction hypothesis is that (**??**) holds if there are $N-1$ or fewer generations. Suppose $L$ has $N$ generations, and let $L_1, \ldots, L_{b_0}$ denote the lattice trees resulting from deleting from $L$ all bonds incident on the origin. We regard each $L_a$ as rooted at the neighbour of the origin in the corresponding deleted bond. It suffices to show that $\nu(L) = b_0! \prod_{a=1}^{b_0} \nu(L_a)$, since each $L_a$ has fewer than $N$ generations.

To prove this, we note that each mean-field configuration $(T,\varphi)$ with $\varphi(T) = L$ induces a set of $(T_a, \varphi_a)$ such that $\varphi_a(T_a) = L_a$. This correspondence is $b_0!$ to 1, since $(T,\varphi)$ is determined by the set of $(T_a, \varphi_a)$, up to permutation of the branches of $T$ at its root. This proves $\nu(L) = b_0! \prod_{a=1}^{b_0} \nu(L_a)$. □

For a general offspring distribution, the above proof gives

$$\sum_{(T,\varphi):\varphi(T)=L} \mathbb{Q}_n^\infty(T,\varphi) = \frac{\prod_{x \in L} p_{b_x} b_x!}{\sum_{L \in \mathcal{L}_n} \prod_{x \in L} p_{b_x} b_x!}, \tag{5.7}$$

which associates degree-dependent weights to each vertex in a lattice tree.

# Acknowledgements

One of us (G.S.) thanks Eric Derbez for early conversations on the feasibility of a combinatorial interpretation of the self-consistent mean-field model introduced in [?]. The authors thank Mike Sachs for focusing our attention on the question of existence of such a combinatorial interpretation. We also thank both Mike Sachs and Jeong Han Kim for valuable conversations. The work of two of us (R.v.d.H. and G.S.) was supported in part by NSERC and was carried out in part during a visit to Microsoft Research in 1998.